\documentclass[11pt]{article}
\usepackage{mathrsfs}
\usepackage{nicefrac}
\usepackage[T1]{fontenc}
\usepackage{latexsym,amssymb,amsmath,amsfonts,amsthm}
\usepackage{graphicx}
\usepackage{caption}
\usepackage[utf8]{inputenc}
\usepackage{authblk}
\usepackage{bm}
\begin{document}

\centerline{\bf\Large New perturbation bounds for the spectrum of a normal matrix}

\medskip
\bigskip

\centerline{Xuefeng Xu$^{{\rm a},{\rm b},\ast}$, Chen-Song Zhang$^{{\rm a},{\rm b}}$}

\medskip

\centerline{\small $^{\rm a}$ Institute of Computational Mathematics and Scientific/Engineering Computing,}
\centerline{\small Academy of Mathematics and Systems Science, Chinese Academy of Sciences, Beijing 100190, China}
\centerline{\small $^{\rm b}$ School of Mathematical Sciences, University of Chinese Academy of Sciences, Beijing 100049, China}
\centerline{\footnotesize $^{\ast}$Corresponding author. E-mail addresses: xuxuefeng@lsec.cc.ac.cn (X. Xu), zhangcs@lsec.cc.ac.cn (C.-S. Zhang).}

\medskip
\bigskip

\noindent{\bf \large Abstract}

Let $A\in\mathbb{C}^{n\times n}$ and $\widetilde{A}\in\mathbb{C}^{n\times n}$ be two normal matrices with spectra $\{\lambda_{i}\}_{i=1}^{n}$ and $\{\widetilde{\lambda}_{i}\}_{i=1}^{n}$, respectively. The celebrated Hoffman--Wielandt theorem states that there exists a permutation $\pi$ of $\{1,\ldots,n\}$ such that $\left(\sum_{i=1}^{n}\big|\widetilde{\lambda}_{\pi(i)}-\lambda_{i}\big|^{2}\right)^{1\over 2}$ is no larger than the Frobenius norm of  $\widetilde{A}-A$. However, if either $A$ or $\widetilde{A}$ is non-normal, this result does not hold in general. In this paper, we present several novel upper bounds for $\left(\sum_{i=1}^{n}\big|\widetilde{\lambda}_{\pi(i)}-\lambda_{i}\big|^{2}\right)^{1\over 2}$, provided that $A$ is normal and $\widetilde{A}$ is arbitrary. Some of these estimates involving the ``departure from normality'' of $\widetilde{A}$ have generalized the Hoffman--Wielandt theorem. Furthermore, we give new perturbation bounds for the spectrum of a Hermitian matrix.

\medskip


\medskip


\medskip
\bigskip
\bigskip

\noindent{\bf \large 1. Introduction}

\medskip

Let $\mathbb{C}^{m\times n}$ and $\mathscr{U}_{n}$ be the set of all $m\times n$ complex matrices and the set of all $n\times n$ unitary matrices, respectively. For any $X\in\mathbb{C}^{m\times n}$, $X^{\ast}$, ${\rm rank}(X)$, $\|X\|_{2}$, and $\|X\|_{F}$ denote the conjugate transpose, the rank, the spectral norm, and the Frobenius norm of $X$, respectively. For any $M\in\mathbb{C}^{n\times n}$, its diagonal part, strictly lower triangular part, and strictly upper triangular part are denoted by $\mathcal{D}(M)$, $\mathcal{L}(M)$, and $\mathcal{U}(M)$, respectively. The trace of $M$ is denoted by ${\rm tr}(M)$. The set $\mathscr{U}_{n}(M)$ is defined as
\begin{displaymath}
\mathscr{U}_{n}(M):=\big\{U\in\mathscr{U}_{n}: U^{\ast}MU \ \text{is upper triangular}\big\}.
\end{displaymath}
The symbol $\kappa_{2}(M)$ stands for the spectral condition number of a nonsingular matrix $M$, namely, $\kappa_{2}(M)=\|M^{-1}\|_{2}\|M\|_{2}$. 

For $M=(m_{ij})\in\mathbb{C}^{n\times n}$, we define $W_{L}(M)$ and $W_{U}(M)$ as follows:
\begin{align*}
W_{L}(M):&=\max\big\{i-j: m_{ij}\neq0, i>j\big\},\tag{1.1a}\\
W_{U}(M):&=\max\big\{j-i: m_{ij}\neq0, i<j\big\}.\tag{1.1b}
\end{align*}
In particular, if $\mathcal{L}(M)=0$ (resp., $\mathcal{U}(M)=0$), we set $W_{L}(M)=0$ (resp., $W_{U}(M)=0$).
Another quantity $\delta(M)$ is defined as
\begin{align*}
\delta(M):=\bigg(\|M\|_{F}^{2}-\frac{1}{n}|{\rm tr}(M)|^{2}\bigg)^{1\over 2},\tag{1.2}
\end{align*}
which is well-defined because $\|M\|_{F}^{2}\geq\frac{1}{n}|{\rm tr}(M)|^{2}$ for all $M\in\mathbb{C}^{n\times n}$. Obviously, $\delta(M)\leq\|M\|_{F}$, and $\delta(M)=\|M\|_{F}$ if and only if ${\rm tr}(M)=0$. For simplicity, for a permutation $\pi$ of $\{1,\ldots,n\}$, we define
\begin{align*}
\mathbb{D}_{2}:=\bigg(\sum_{i=1}^{n}\big|\widetilde{\lambda}_{\pi(i)}-\lambda_{i}\big|^{2}\bigg)^{1\over 2},\tag{1.3}
\end{align*}
which characterizes the ``distance'' between $\{\lambda_{i}\}_{i=1}^{n}$ and $\{\widetilde{\lambda}_{i}\}_{i=1}^{n}$ with respect to $\ell^{2}$-norm. 

Recall that $A\in\mathbb{C}^{n\times n}$ is \emph{normal} if $A$ commutes with its conjugate transpose, i.e., $AA^{\ast}=A^{\ast}A$. In particular, if $A=A^{\ast}$, then $A\in\mathbb{C}^{n\times n}$ is called a \emph{Hermitian} matrix. Assume that $A\in\mathbb{C}^{n\times n}$ and $\widetilde{A}\in\mathbb{C}^{n\times n}$ are normal matrices with spectra $\{\lambda_{i}\}_{i=1}^{n}$ and $\{\widetilde{\lambda}_{i}\}_{i=1}^{n}$, respectively. Let $E=\widetilde{A}-A$. In 1953, Hoffman and Wielandt~[1] proved that there is a permutation $\pi$ of $\{1,\ldots,n\}$ such that
\begin{displaymath}
\mathbb{D}_{2}\leq\|E\|_{F},
\end{displaymath}
which is the well-known \emph{Hoffman--Wielandt theorem}. This theorem reveals that there is a strong global stability to the set of eigenvalues of a normal matrix. Unfortunately, the inequality may fail when $\widetilde{A}$ is non-normal. For example, if
\begin{displaymath}
A=\begin{pmatrix}
0 & 0 \\
0 & 3
\end{pmatrix} \quad \text{and} \quad \widetilde{A}=\begin{pmatrix}
-1 & -1 \\
1 & 1
\end{pmatrix},
\end{displaymath}
then $A$ is a normal matrix with spectrum $\{3,0\}$ and $\widetilde{A}$ is a non-normal matrix with spectrum $\{0,0\}$. For any permutation $\pi$ of $\{1,2\}$, we have that $\mathbb{D}_{2}=3>\sqrt{7}=\|E\|_{F}$.

Due to the limitation of the Hoffman--Wielandt theorem, many authors have developed analogous results. If $A\in\mathbb{C}^{n\times n}$ is normal and $\widetilde{A}\in\mathbb{C}^{n\times n}$ is non-normal, Sun~[2] demonstrated that there exists a permutation $\pi$ of $\{1,\ldots,n\}$ such that
\begin{align*}
\mathbb{D}_{2}\leq\sqrt{n}\|E\|_{F},\tag{1.4}
\end{align*}
which can be applied to characterize the variation of the spectrum of an arbitrary matrix [3]. If $A\in\mathbb{C}^{n\times n}$ is normal and $\widetilde{A}\in\mathbb{C}^{n\times n}$ can be diagonalized by $\widetilde{X}$, i.e., $\widetilde{A}=\widetilde{X}\widetilde{\Lambda}\widetilde{X}^{-1}$ ($\widetilde{\Lambda}$ is diagonal), Sun~[4] proved that there exists a permutation $\pi$ of $\{1,\ldots,n\}$ such that
\begin{displaymath}
\mathbb{D}_{2}\leq\kappa_{2}(\widetilde{X})\|E\|_{F}.
\end{displaymath}
If both $A\in\mathbb{C}^{n\times n}$ and $\widetilde{A}\in\mathbb{C}^{n\times n}$ are diagonalizable, i.e., there are nonsingular matrices $S\in\mathbb{C}^{n\times n}$ and $T\in\mathbb{C}^{n\times n}$ such that $S^{-1}AS={\rm diag}\big(\lambda_{1},\ldots,\lambda_{n}\big)$ and $T^{-1}\widetilde{A}T={\rm diag}\big(\widetilde{\lambda}_{1},\ldots,\widetilde{\lambda}_{n}\big)$, Zhang~[5] showed that there exists a permutation $\pi$ of $\{1,\ldots,n\}$ such that
\begin{displaymath}
\mathbb{D}_{2}\leq\kappa_{2}(S)\kappa_{2}(T)\|E\|_{F}.
\end{displaymath}
Assume that $A\in\mathbb{C}^{n\times n}$ is normal, $\widetilde{A}\in\mathbb{C}^{n\times n}$ is arbitrary, and $\widetilde{U}_{1}\in\mathscr{U}_{n}$ such that 
\begin{displaymath}
\widetilde{U}_{1}^{\ast}\widetilde{A}\widetilde{U}_{1}={\rm diag}\big(\widetilde{A}_{1},\ldots,\widetilde{A}_{s}\big)
\end{displaymath}
for some positive integer $s$, where each $\widetilde{A}_{i}\in\mathbb{C}^{n_{i}\times n_{i}}$ is upper triangular and $\sum_{i=1}^{s}n_{i}=n$. Li and Sun~[6] proved that there exists a permutation $\pi$ of $\{1,\ldots,n\}$ such that
\begin{align*}
\mathbb{D}_{2}\leq\sqrt{n-s+1}\|E\|_{F},\tag{1.5}
\end{align*}
which has improved (1.4). In particular, if $\widetilde{A}$ is normal (hence $s=n$), then (1.5) reduces to the Hoffman--Wielandt theorem. As a special case of normal matrices, if $A\in\mathbb{C}^{n\times n}$ is Hermitian and $\widetilde{A}\in\mathbb{C}^{n\times n}$ is non-normal, then by Kahan's result~[7], there exists a permutation $\pi$ of $\{1,\ldots,n\}$ such that
\begin{align*}
\mathbb{D}_{2}\leq\sqrt{2}\|E\|_{F}.\tag{1.6}
\end{align*}

It is easy to see that the upper bounds in (1.4) and (1.6) only depend on the distance $\|E\|_{F}$. In other words, the bounds do not change no matter how close  $\widetilde{A}\widetilde{A}^{\ast}$ is to $\widetilde{A}^{\ast}\widetilde{A}$. By the Schur’s theorem, there exists a $\widetilde{U}\in\mathscr{U}_{n}$ such that
\begin{displaymath} \widetilde{A}=\widetilde{U}\big(\widetilde{\Lambda}+\Delta\big)\widetilde{U}^{\ast},
\end{displaymath}
where $\widetilde{\Lambda}={\rm diag}\big(\widetilde{\lambda}_{1},\ldots,\widetilde{\lambda}_{n}\big)$ and $\Delta$ is strictly upper triangular. Hence,
\begin{displaymath}
\|\Delta\|_{F}=\bigg(\|\widetilde{A}\|_{F}^{2}-\sum_{i=1}^{n}\big|\widetilde{\lambda}_{i}\big|^{2}\bigg)^{1\over 2},
\end{displaymath}
which can be considered as a quantitative measure of the non-normality of $\widetilde{A}$. Indeed, $\|\Delta\|_{F}$ is referred to as the \emph{departure from normality} (with respect to $\|\cdot\|_{F}$) of $\widetilde{A}$ [8]. Taking the quantity into account, Sun~[9] established that there exists a permutation $\pi$ of $\{1,\ldots,n\}$ such that
\begin{align*}
\mathbb{D}_{2}\leq\sqrt{\|E\|_{F}^{2}+2\min\big\{\|A\|_{F},\sqrt{n-1}\|A\|_{2}\big\}\|\Delta\|_{F}-\|\Delta\|_{F}^{2}},\tag{1.7}
\end{align*}
provided that $A\in\mathbb{C}^{n\times n}$ is normal and $\widetilde{A}\in\mathbb{C}^{n\times n}$ is arbitrary. Recently, Li and Vong~[10] studied the variation of the spectrum of a Hermitian matrix and obtained that there is a permutation $\pi$ of $\{1,\ldots,n\}$ such that
\begin{align*}
\mathbb{D}_{2}&\leq\sqrt{\|E\|_{F}^{2}+\sqrt{2}\|E\|_{F}\|\Delta\|_{F}},\tag{1.8}\\
\mathbb{D}_{2}&\leq\sqrt{\|E\|_{F}^{2}+2\|E\|_{F}\|\Delta\|_{F}-\|\Delta\|_{F}^{2}},\tag{1.9}
\end{align*}
provided that $A\in\mathbb{C}^{n\times n}$ is Hermitian and $\widetilde{A}\in\mathbb{C}^{n\times n}$ is arbitrary. 

In this paper, we establish some novel estimates for $\mathbb{D}_{2}$ (see Theorems 3.6 and 3.10 below), provided that $A\in\mathbb{C}^{n\times n}$ is normal and $\widetilde{A}\in\mathbb{C}^{n\times n}$ is arbitrary. The main results include:
\begin{align*}
\mathbb{D}_{2}&\leq\sqrt{\|E\|_{F}^{2}+(n-1)\delta(E)^{2}},\tag{1.10}\\
\mathbb{D}_{2}&\leq\sqrt{\|E\|_{F}^{2}+(n-s)\delta(E)^{2}},\tag{1.11}\\
\mathbb{D}_{2}&\leq\sqrt{\|E\|_{F}^{2}+2\sqrt{\frac{n-1}{n}}\delta(A)\|\Delta\|_{F}-\|\Delta\|_{F}^{2}},\tag{1.12}
\end{align*}
where $s$ is defined as in (1.5). Note that, if ${\rm tr}(E)\neq0$ (hence $\delta(E)<\|E\|_{F}$), our estimates (1.10) and (1.11) are sharper than (1.4) and (1.5), respectively. Since $\|A\|_{F}^{2}\leq n\|A\|_{2}^{2}$, we have
\begin{displaymath}
\sqrt{\frac{n-1}{n}}\delta(A)\leq\min\big\{\|A\|_{F},\sqrt{n-1}\|A\|_{2}\big\},
\end{displaymath}
which implies that the upper bound in (1.12) is smaller than that in (1.7). Moreover, we establish some new perturbation bounds for the spectrum of a Hermitian matrix (see Theorem~4.2 below), which contain:
\begin{align*}
\mathbb{D}_{2}&\leq\sqrt{\|E\|_{F}^{2}+\delta(E)^{2}},\tag{1.13}\\
\mathbb{D}_{2}&\leq\sqrt{\|E\|_{F}^{2}+\sqrt{2}\delta(E)\|\Delta\|_{F}},\tag{1.14}\\
\mathbb{D}_{2}&\leq\sqrt{\|E\|_{F}^{2}+2\delta(E)\|\Delta\|_{F}-\|\Delta\|_{F}^{2}}.\tag{1.15}
\end{align*}
Similarly, if ${\rm tr}(E)=0$, then our results (1.13)--(1.15) coincide with (1.6), (1.8), and (1.9), respectively; otherwise, our upper bounds are smaller.  

The rest of this paper is organized as follows. In Section 2, we introduce some relations between $\|\mathcal{U}(A)\|_{F}$ and $\|\mathcal{L}(A)\|_{F}$, provided that $A$ is a normal matrix. In Section 3, we establish some novel perturbation bounds for the spectrum of a normal matrix.  In Section 4, we present new perturbation bounds for the spectrum of a Hermitian matrix. Finally, some conclusions are given in Section 5. 

\medskip
\bigskip

\noindent{\bf \large 2. Preliminaries}

\medskip

In this section, we introduce several useful properties of normal matrices. The following lemma gives an identity on the entries of a normal matrix (see, e.g., [2, Lemma 2.1]).

\medskip

\noindent{\bf Lemma 2.1.} \emph{Let $A=(a_{ij})\in\mathbb{C}^{n\times n}$ be normal. Then}
\begin{displaymath}
\sum_{i=1}^{n-1}\sum_{j=i+1}^{n}(j-i)|a_{ij}|^{2}=\sum_{j=1}^{n-1}\sum_{i=j+1}^{n}(i-j)|a_{ij}|^{2}.
\end{displaymath}

\medskip

Using Lemma 2.1, we can obtain the following relations between $\|\mathcal{U}(A)\|_{F}$ and $\|\mathcal{L}(A)\|_{F}$.

\medskip

\noindent{\bf Lemma 2.2.} \emph{Let $A=(a_{ij})\in\mathbb{C}^{n\times n}$ be normal. Then}
\begin{align*}
\|\mathcal{U}(A)\|_{F}&\leq\sqrt{W_{L}(A)}\|\mathcal{L}(A)\|_{F},\tag{2.1a}\\ \|\mathcal{L}(A)\|_{F}&\leq\sqrt{W_{U}(A)}\|\mathcal{U}(A)\|_{F}.\tag{2.1b}
\end{align*}

\noindent{\bf Proof.} According to the definition of $\mathcal{U}(\cdot)$, it follows that
\begin{displaymath}
\|\mathcal{U}(A)\|_{F}^{2}=\sum_{i=1}^{n-1}\sum_{j=i+1}^{n}|a_{ij}|^{2}\leq\sum_{i=1}^{n-1}\sum_{j=i+1}^{n}(j-i)|a_{ij}|^{2}.
\end{displaymath}
By Lemma 2.1, we have
\begin{displaymath}
\|\mathcal{U}(A)\|_{F}^{2}\leq\sum_{j=1}^{n-1}\sum_{i=j+1}^{n}(i-j)|a_{ij}|^{2}\leq W_{L}(A)\sum_{j=1}^{n-1}\sum_{i=j+1}^{n}|a_{ij}|^{2}=W_{L}(A)\|\mathcal{L}(A)\|_{F}^{2},
\end{displaymath}
which leads to (2.1a). Analogously, we can prove the second inequality. \qed

\medskip

\noindent{\bf Remark 2.3.} For any $A\in\mathbb{C}^{n\times n}$, it is clear that
\begin{displaymath}
W_{L}(A)\leq n-1 \quad \text{and} \quad W_{U}(A)\leq n-1.
\end{displaymath}
Hence, from Lemma 2.2, we obtain
\begin{align*}
\|\mathcal{U}(A)\|_{F}&\leq\sqrt{n-1}\|\mathcal{L}(A)\|_{F},\tag{2.2a}\\
\|\mathcal{L}(A)\|_{F}&\leq\sqrt{n-1}\|\mathcal{U}(A)\|_{F},\tag{2.2b}
\end{align*}
which are the inequalities stated in [2, Lemma 3.1].

\medskip

The following lemma presents a modified version of (2.2a) and (2.2b). For a detailed proof, we refer the interested reader to [6, Lemma 2.2].

\medskip

\noindent{\bf Lemma 2.4.} \emph{Let $A\in\mathbb{C}^{n\times n}$ be normal. Then for any $i\in\{1,\dots,n\}$,}
\begin{displaymath}
\|(\mathcal{U}(A))_{(i)}\|_{F}\leq\|\mathcal{L}(A)\|_{F} \quad \emph{and} \quad
\|(\mathcal{L}(A))_{(i)}\|_{F}\leq\|\mathcal{U}(A)\|_{F},
\end{displaymath}
\emph{where $(\cdot)_{(i)}$ denotes the $i$-th row of a matrix.}

\medskip
\bigskip

\noindent{\bf \large 3. Perturbation bounds for the spectrum of a normal matrix}

\medskip

In this section, we present several novel perturbation bounds for the spectrum of a normal matrix. Some of our estimates have improved the existing results in [2, 6, 9].
 
\bigskip

{\bf 3.1. Two useful lemmas.} We first prove a simple but important lemma, which plays a key role in our further analysis.

\medskip

\noindent{\bf Lemma 3.1.} \emph{Let $M=(m_{ij})\in\mathbb{C}^{n\times n}$. Then}
\begin{displaymath}
\|\mathcal{L}(M)\|_{F}^{2}+\|\mathcal{U}(M)\|_{F}^{2}\leq\delta(M)^{2}.
\end{displaymath}

\noindent{\bf Proof.} According to the fact that
\begin{displaymath}
\|M\|_{F}^{2}=\|\mathcal{D}(M)\|_{F}^{2}+\|\mathcal{L}(M)\|_{F}^{2}+\|\mathcal{U}(M)\|_{F}^{2},
\end{displaymath}
using the Cauchy--Schwarz’s inequality, we immediately obtain
\begin{displaymath}
\|\mathcal{L}(M)\|_{F}^{2}+\|\mathcal{U}(M)\|_{F}^{2}=\|M\|_{F}^{2}-\sum_{i=1}^{n}|m_{ii}|^{2}\leq\|M\|_{F}^{2}-\frac{1}{n}\bigg(\sum_{i=1}^{n}|m_{ii}|\bigg)^{2}.
\end{displaymath}
Note that
\begin{displaymath} \sum_{i=1}^{n}|m_{ii}|\geq\bigg|\sum_{i=1}^{n}m_{ii}\bigg|=|{\rm tr}(M)|.
\end{displaymath}
We then get the desired inequality. \qed

\medskip

For any $M\in\mathbb{C}^{n\times n}$, one may give other upper bounds for $\|\mathcal{L}(M)\|_{F}^{2}+\|\mathcal{U}(M)\|_{F}^{2}$. It is well-known that the \emph{Hadamard product} of $A=(a_{ij})\in\mathbb{C}^{m\times n}$ and $B=(b_{ij})\in\mathbb{C}^{m\times n}$ is defined as $A\circ B=(a_{ij}b_{ij})\in\mathbb{C}^{m\times n}$ (see, e.g., [11, Definition 5.0.1]). According to the proof of Lemma 3.1, we have
\begin{displaymath}
\|\mathcal{L}(M)\|_{F}^{2}+\|\mathcal{U}(M)\|_{F}^{2}=\|M\|_{F}^{2}-\sum_{i=1}^{n}|m_{ii}^{2}|\leq\|M\|_{F}^{2}-\bigg|\sum_{i=1}^{n}m_{ii}^{2}\bigg|,
\end{displaymath}
which yields 
\begin{displaymath}
\|\mathcal{L}(M)\|_{F}^{2}+\|\mathcal{U}(M)\|_{F}^{2}\leq\|M\|_{F}^{2}-|{\rm tr}(M\circ M)|.
\end{displaymath}
For any $M=(m_{ij})\in\mathbb{C}^{n\times n}$, the \emph{entry-wise absolute value} of $M$ is defined as $|M|=(|m_{ij}|)\in\mathbb{R}_{+}^{n\times n}$ (see, e.g., [11, p. 124]), where $\mathbb{R}_{+}^{n\times n}$ denotes the set of all $n\times n$ non-negative matrices. We can also show that
\begin{align*}
\|\mathcal{L}(M)\|_{F}^{2}+\|\mathcal{U}(M)\|_{F}^{2}&=\|M\|_{F}^{2}-{\rm tr}(|M|\circ|M|),\\
\|\mathcal{L}(M)\|_{F}^{2}+\|\mathcal{U}(M)\|_{F}^{2}&\leq\|M\|_{F}^{2}-\frac{1}{n}\left({\rm tr}(|M|)\right)^{2}.
\end{align*}

In view of the above relations, we now define
\begin{align*}
\phi_{1}(M):&=\|M\|_{F}^{2}-|{\rm tr}(M\circ M)|,\\ \phi_{2}(M):&=\|M\|_{F}^{2}-{\rm tr}(|M|\circ|M|),\\ \phi_{3}(M):&=\|M\|_{F}^{2}-\frac{1}{n}\left({\rm tr}(|M|)\right)^{2}.
\end{align*}
It is clear that $\delta(M)$ is invariant under a unitary similarity transformation, i.e.,
\begin{displaymath} \delta(U^{\ast}MU)=\delta(M)
\end{displaymath}
for all $U\in\mathscr{U}_{n}$. However, $\phi_{j}(\cdot) \ (j=1,2,3)$ may change under a unitary similarity transformation. For instance,
\begin{displaymath}
M=\begin{pmatrix}
1+i & 0\\
0 & 2
\end{pmatrix} \quad \text{and} \quad U=\frac{1}{\sqrt{2}}\begin{pmatrix}
1 & i\\
i & 1
\end{pmatrix},
\end{displaymath}
where $i=\sqrt{-1}$. Then we have
\begin{displaymath}
U^{\ast}MU=\frac{1}{2}\begin{pmatrix}
3+i & -1-i\\
1+i & 3+i
\end{pmatrix} \quad \text{and} \quad |U^{\ast}MU|=\frac{1}{2}\begin{pmatrix}
\sqrt{10} & \sqrt{2}\\
\sqrt{2} & \sqrt{10}
\end{pmatrix}.
\end{displaymath}
Straightforward calculations yield
\begin{align*}
\phi_{1}(M)&=6-2\sqrt{5}>1=\phi_{1}(U^{\ast}MU),\\
\phi_{2}(M)&=0<1=\phi_{2}(U^{\ast}MU),\\
\phi_{3}(M)&=3-2\sqrt{2}<1=\phi_{3}(U^{\ast}MU).
\end{align*}
This example also illustrates that $\phi_{j}(M) \ (j=1,2,3)$ may alter under a unitary similarity transformation, even if $M$ is a normal matrix.

\medskip

Using Lemmas 2.2 and 3.1, we can get the following estimates for $\|\mathcal{U}(A)\|_{F}$ and $\|\mathcal{L}(A)\|_{F}$.

\medskip

\noindent{\bf Lemma 3.2.} \emph{Let $A\in\mathbb{C}^{n\times n}$ be normal. Then}
\begin{displaymath}
\|\mathcal{U}(A)\|_{F}\leq \sqrt{\frac{W_{L}(A)}{1+W_{L}(A)}}\delta(A) \quad \emph{and} \quad
\|\mathcal{L}(A)\|_{F}\leq \sqrt{\frac{W_{U}(A)}{1+W_{U}(A)}}\delta(A).
\end{displaymath}

\noindent{\bf Proof.} By (2.1a), we have
\begin{displaymath}
\left(1+W_{L}(A)\right)\|\mathcal{U}(A)\|_{F}^{2}\leq W_{L}(A)\left(\|\mathcal{L}(A)\|_{F}^{2}+\|\mathcal{U}(A)\|_{F}^{2}\right).
\end{displaymath}
Using Lemma 3.1, we obtain
\begin{displaymath}
\left(1+W_{L}(A)\right)\|\mathcal{U}(A)\|_{F}^{2}\leq W_{L}(A)\delta(A)^{2},
\end{displaymath}
which yields
\begin{displaymath}
\|\mathcal{U}(A)\|_{F}\leq\sqrt{ \frac{W_{L}(A)}{1+W_{L}(A)}}\delta(A).
\end{displaymath}
Similarly, by (2.1b) and Lemma 3.1, we have
\begin{displaymath}
\|\mathcal{L}(A)\|_{F}\leq\sqrt{ \frac{W_{U}(A)}{1+W_{U}(A)}}\delta(A).
\end{displaymath}
This completes the proof. \qed

\medskip

From Remark 2.3, we can see that the inequalities in Lemma 3.2 still hold when $W_{L}(A)$ and $W_{U}(A)$ are replaced by $n-1$. This observation yields the following corollary.

\medskip

\noindent{\bf Corollary 3.3.} \emph{Let $A\in\mathbb{C}^{n\times n}$ be normal. Then}
\begin{displaymath}
\max\bigg\{\sup_{U\in\mathscr{U}_{n}}\|\mathcal{U}(U^{\ast}AU)\|_{F},\sup_{U\in\mathscr{U}_{n}}\|\mathcal{L}(U^{\ast}AU)\|_{F}\bigg\}\leq \sqrt{\frac{n-1}{n}}\delta(A).
\end{displaymath}
\emph{In particular, we have}
\begin{displaymath}
\max\big\{\|\mathcal{U}(A)\|_{F},\|\mathcal{L}(A)\|_{F}\big\}\leq \sqrt{\frac{n-1}{n}}\delta(A).
\end{displaymath}

\bigskip

{\bf 3.2. The estimates based on Schur’s decomposition.} Let $A\in\mathbb{C}^{n\times n}$ be a normal matrix with spectrum $\{\lambda_{i}\}_{i=1}^{n}$. Assume that $\widetilde{A}=A+E$ has the spectrum $\{\widetilde{\lambda}_{i}\}_{i=1}^{n}$, where $E\in\mathbb{C}^{n\times n}$ is an arbitrary perturbation. By the Schur's theorem, there exists a $\widetilde{U}\in\mathscr{U}_{n}$ such that
\begin{displaymath}
\widetilde{A}=\widetilde{U}\big(\widetilde{\Lambda}+\Delta\big)\widetilde{U}^{\ast},
\end{displaymath}
where $\widetilde{\Lambda}$ is diagonal and $\Delta$ is strictly upper triangular. From $\widetilde{A}=A+E$, we have
\begin{displaymath}
\widetilde{U}^{\ast}A\widetilde{U}+\widetilde{U}^{\ast}E\widetilde{U}=\widetilde{\Lambda}+\Delta,
\end{displaymath} 
which means
\begin{align*}
\mathcal{L}(\widetilde{U}^{\ast}A\widetilde{U})+\mathcal{L}(\widetilde{U}^{\ast}E\widetilde{U})&=0,\tag{3.1a}\\
\mathcal{U}(\widetilde{U}^{\ast}A\widetilde{U})+\mathcal{U}(\widetilde{U}^{\ast}E\widetilde{U})&=\Delta.\tag{3.1b}
\end{align*}
Since both $\widetilde{\Lambda}$ and $\widetilde{U}^{\ast}A\widetilde{U}$ are normal, by the Hoffman--Wielandt theorem, we obtain that there exists a permutation $\pi$ of $\{1,\ldots,n\}$ such that
\begin{displaymath}
\mathbb{D}_{2}\leq\|\widetilde{\Lambda}-\widetilde{U}^{\ast}A\widetilde{U}\|_{F}=\|\widetilde{U}^{\ast}E\widetilde{U}-\Delta\|_{F}.
\end{displaymath}
It follows from (3.1b) that $\widetilde{U}^{\ast}E\widetilde{U}-\Delta$ can be written as
\begin{displaymath}
\widetilde{U}^{\ast}E\widetilde{U}-\Delta=\mathcal{D}(\widetilde{U}^{\ast}E\widetilde{U})+\mathcal{L}(\widetilde{U}^{\ast}E\widetilde{U})-\mathcal{U}(\widetilde{U}^{\ast}A\widetilde{U}).
\end{displaymath}
Then,
\begin{align*}
\|\widetilde{U}^{\ast}E\widetilde{U}-\Delta\|_{F}^{2}&=\|\mathcal{D}(\widetilde{U}^{\ast}E\widetilde{U})\|_{F}^{2}+\|\mathcal{L}(\widetilde{U}^{\ast}E\widetilde{U})\|_{F}^{2}+\|\mathcal{U}(\widetilde{U}^{\ast}A\widetilde{U})\|_{F}^{2}\\
&=\|E\|_{F}^{2}+\|\mathcal{U}(\widetilde{U}^{\ast}A\widetilde{U})\|_{F}^{2}-\|\mathcal{U}(\widetilde{U}^{\ast}E\widetilde{U})\|_{F}^{2}.
\end{align*}
Hence, we obtain that there exists a permutation $\pi$ of $\{1,\ldots,n\}$ such that
\begin{align*}
\mathbb{D}_{2}\leq\sqrt{\|E\|_{F}^{2}+\|\mathcal{U}(\widetilde{U}^{\ast}A\widetilde{U})\|_{F}^{2}-\|\mathcal{U}(\widetilde{U}^{\ast}E\widetilde{U})\|_{F}^{2}}.\tag{3.2}
\end{align*}

In order to derive the upper bounds for $\mathbb{D}_{2}$, we need to estimate $\|\mathcal{U}(\widetilde{U}^{\ast}A\widetilde{U})\|_{F}^{2}-\|\mathcal{U}(\widetilde{U}^{\ast}E\widetilde{U})\|_{F}^{2}$. Based on the different estimates for $\|\mathcal{U}(\widetilde{U}^{\ast}A\widetilde{U})\|_{F}^{2}-\|\mathcal{U}(\widetilde{U}^{\ast}E\widetilde{U})\|_{F}^{2}$, we can obtain the following lemma.

\medskip

\noindent{\bf Lemma 3.4.} \emph{Let $A\in\mathbb{C}^{n\times n}$ be a normal matrix with spectrum $\{\lambda_{i}\}_{i=1}^{n}$, and let $\widetilde{A}=A+E$ with spectrum $\{\widetilde{\lambda}_{i}\}_{i=1}^{n}$, where $E\in\mathbb{C}^{n\times n}$ is an arbitrary perturbation. Then there exists a permutation $\pi$ of $\{1,\ldots,n\}$ such that}
\begin{align*}
\mathbb{D}_{2}&\leq\sqrt{\|E\|_{F}^{2}+W_{L}(\widetilde{U}^{\ast}E\widetilde{U})\delta(E)^{2}},\tag{3.3a}\\
\mathbb{D}_{2}&\leq\sqrt{\|E\|_{F}^{2}+\sqrt{1+W_{L}(\widetilde{U}^{\ast}E\widetilde{U})}\delta(E)\|\Delta\|_{F}},\tag{3.3b}\\
\mathbb{D}_{2}&\leq\sqrt{\|E\|_{F}^{2}+2\delta(E)\|\Delta\|_{F}+\|\Delta\|_{F}^{2}},\tag{3.3c}\\
\mathbb{D}_{2}&\leq\sqrt{\|E\|_{F}^{2}+2\sqrt{W_{L}(\widetilde{U}^{\ast}E\widetilde{U})}\delta(E)\|\Delta\|_{F}-\|\Delta\|_{F}^{2}},\tag{3.3d}
\end{align*}
\emph{where $\Delta=\mathcal{U}(\widetilde{U}^{\ast}\widetilde{A}\widetilde{U})$ with $\widetilde{U}\in\mathscr{U}_{n}(\widetilde{A})$.}

\medskip

\noindent{\bf Proof.} 
(a) Since $\widetilde{U}^{\ast}A\widetilde{U}$ is normal, by (2.1a), we have
\begin{displaymath}
\|\mathcal{U}(\widetilde{U}^{\ast}A\widetilde{U})\|_{F}^{2}-\|\mathcal{U}(\widetilde{U}^{\ast}E\widetilde{U})\|_{F}^{2}\leq W_{L}(\widetilde{U}^{\ast}A\widetilde{U})\|\mathcal{L}(\widetilde{U}^{\ast}A\widetilde{U})\|_{F}^{2}-\|\mathcal{U}(\widetilde{U}^{\ast}E\widetilde{U})\|_{F}^{2}.
\end{displaymath}
According to (3.1a) and the definition of $W_{L}(\cdot)$, it follows that $W_{L}(\widetilde{U}^{\ast}A\widetilde{U})=W_{L}(\widetilde{U}^{\ast}E\widetilde{U})$. Hence,
\begin{align*}
\|\mathcal{U}(\widetilde{U}^{\ast}A\widetilde{U})\|_{F}^{2}-\|\mathcal{U}(\widetilde{U}^{\ast}E\widetilde{U})\|_{F}^{2}&\leq W_{L}(\widetilde{U}^{\ast}E\widetilde{U})\|\mathcal{L}(\widetilde{U}^{\ast}E\widetilde{U})\|_{F}^{2}-\|\mathcal{U}(\widetilde{U}^{\ast}E\widetilde{U})\|_{F}^{2}\\
&\leq W_{L}(\widetilde{U}^{\ast}E\widetilde{U})\Big(\|\mathcal{L}(\widetilde{U}^{\ast}E\widetilde{U})\|_{F}^{2}+\|\mathcal{U}(\widetilde{U}^{\ast}E\widetilde{U})\|_{F}^{2}\Big)\\
&\leq W_{L}(\widetilde{U}^{\ast}E\widetilde{U})\delta(E)^{2},
\end{align*}
where we have used Lemma 3.1. Then the estimate (3.3a) follows immediately from (3.2).

(b) Based on (2.1a), (3.1a), the Cauchy--Schwarz’s inequality, and Lemma 3.1, we obtain
\begin{align*}
\|\mathcal{U}(\widetilde{U}^{\ast}A\widetilde{U})\|_{F}+\|\mathcal{U}(\widetilde{U}^{\ast}E\widetilde{U})\|_{F}&\leq\sqrt{W_{L}(\widetilde{U}^{\ast}E\widetilde{U})}\|\mathcal{L}(\widetilde{U}^{\ast}E\widetilde{U})\|_{F}+\|\mathcal{U}(\widetilde{U}^{\ast}E\widetilde{U})\|_{F}\\
&\leq\sqrt{1+W_{L}(\widetilde{U}^{\ast}E\widetilde{U})}\sqrt{\|\mathcal{L}(\widetilde{U}^{\ast}E\widetilde{U})\|_{F}^{2}+\|\mathcal{U}(\widetilde{U}^{\ast}E\widetilde{U})\|_{F}^{2}}\\
&\leq\sqrt{1+W_{L}(\widetilde{U}^{\ast}E\widetilde{U})}\delta(E).
\end{align*}
From (3.1b), we have 
\begin{displaymath}
\|\mathcal{U}(\widetilde{U}^{\ast}A\widetilde{U})\|_{F}-\|\mathcal{U}(\widetilde{U}^{\ast}E\widetilde{U})\|_{F}\leq\|\Delta\|_{F}.
\end{displaymath}
Hence, the estimate (3.3b) holds because of (3.2).

(c) By (3.1b), the triangle inequality, and Lemma 3.1, we have
\begin{align*}
\|\mathcal{U}(\widetilde{U}^{\ast}A\widetilde{U})\|_{F}^{2}-\|\mathcal{U}(\widetilde{U}^{\ast}E\widetilde{U})\|_{F}^{2}&=\|\Delta-\mathcal{U}(\widetilde{U}^{\ast}E\widetilde{U})\|_{F}^{2}-\|\mathcal{U}(\widetilde{U}^{\ast}E\widetilde{U})\|_{F}^{2}\\
&\leq\Big(\|\Delta\|_{F}+\|\mathcal{U}(\widetilde{U}^{\ast}E\widetilde{U})\|_{F}\Big)^{2}-\|\mathcal{U}(\widetilde{U}^{\ast}E\widetilde{U})\|_{F}^{2}\\
&=2\|\mathcal{U}(\widetilde{U}^{\ast}E\widetilde{U})\|_{F}\|\Delta\|_{F}+\|\Delta\|_{F}^{2}\\
&\leq 2\sqrt{\|\mathcal{L}(\widetilde{U}^{\ast}E\widetilde{U})\|_{F}^{2}+\|\mathcal{U}(\widetilde{U}^{\ast}E\widetilde{U})\|_{F}^{2}}\|\Delta\|_{F}+\|\Delta\|_{F}^{2}\\
&\leq 2\delta(E)\|\Delta\|_{F}+\|\Delta\|_{F}^{2}.
\end{align*}
An application of (3.2) yields the estimate (3.3c).

(d) Using (3.1b) and the triangle inequality, we obtain
\begin{align*}
\|\mathcal{U}(\widetilde{U}^{\ast}A\widetilde{U})\|_{F}^{2}-\|\mathcal{U}(\widetilde{U}^{\ast}E\widetilde{U})\|_{F}^{2}&=\|\mathcal{U}(\widetilde{U}^{\ast}A\widetilde{U})\|_{F}^{2}-\|\Delta-\mathcal{U}(\widetilde{U}^{\ast}A\widetilde{U})\|_{F}^{2}\\
&\leq\|\mathcal{U}(\widetilde{U}^{\ast}A\widetilde{U})\|_{F}^{2}-\Big(\|\Delta\|_{F}-\|\mathcal{U}(\widetilde{U}^{\ast}A\widetilde{U})\|_{F}\Big)^{2}\\
&=2\|\mathcal{U}(\widetilde{U}^{\ast}A\widetilde{U})\|_{F}\|\Delta\|_{F}-\|\Delta\|_{F}^{2}.
\end{align*}
By (2.1a) and (3.1a), we have
\begin{displaymath}
\|\mathcal{U}(\widetilde{U}^{\ast}A\widetilde{U})\|_{F}^{2}-\|\mathcal{U}(\widetilde{U}^{\ast}E\widetilde{U})\|_{F}^{2}\leq 2\sqrt{W_{L}(\widetilde{U}^{\ast}E\widetilde{U})}\|\mathcal{L}(\widetilde{U}^{\ast}E\widetilde{U})\|_{F}\|\Delta\|_{F}-\|\Delta\|_{F}^{2}.
\end{displaymath}
Because $\|\mathcal{L}(\widetilde{U}^{\ast}E\widetilde{U})\|_{F}\leq\sqrt{\|\mathcal{L}(\widetilde{U}^{\ast}E\widetilde{U})\|_{F}^{2}+\|\mathcal{U}(\widetilde{U}^{\ast}E\widetilde{U})\|_{F}^{2}}\leq\delta(E)$, we arrive at
\begin{displaymath}
\|\mathcal{U}(\widetilde{U}^{\ast}A\widetilde{U})\|_{F}^{2}-\|\mathcal{U}(\widetilde{U}^{\ast}E\widetilde{U})\|_{F}^{2}\leq 2\sqrt{W_{L}(\widetilde{U}^{\ast}E\widetilde{U})}\delta(E)\|\Delta\|_{F}-\|\Delta\|_{F}^{2}.
\end{displaymath}
It follows from (3.2) that the estimate (3.3d) holds. \qed

\medskip

We next give another two upper
bounds for $\mathbb{D}_{2}$, which are related to the original matrix $A$.

\medskip

\noindent{\bf Lemma 3.5.} \emph{Under the conditions of Lemma 3.4. The permutation $\pi$ in Lemma 3.4 also~satisfies}
\begin{align*}
\mathbb{D}_{2}&\leq\sqrt{\|E\|_{F}^{2}+\frac{W_{L}(\widetilde{U}^{\ast}E\widetilde{U})}{1+W_{L}(\widetilde{U}^{\ast}E\widetilde{U})}\delta(A)^{2}},\tag{3.4a}\\
\mathbb{D}_{2}&\leq\sqrt{\|E\|_{F}^{2}+2\sqrt{\frac{W_{L}(\widetilde{U}^{\ast}E\widetilde{U})}{1+W_{L}(\widetilde{U}^{\ast}E\widetilde{U})}}\delta(A)\|\Delta\|_{F}-\|\Delta\|_{F}^{2}}.\tag{3.4b}
\end{align*}

\noindent{\bf Proof.} (a) Using Lemma 3.2, we have
\begin{displaymath}
\|\mathcal{U}(\widetilde{U}^{\ast}A\widetilde{U})\|_{F}^{2}\leq \frac{W_{L}(\widetilde{U}^{\ast}A\widetilde{U})}{1+W_{L}(\widetilde{U}^{\ast}A\widetilde{U})}\delta(A)^{2},
\end{displaymath}
which leads to
\begin{displaymath}
\|\mathcal{U}(\widetilde{U}^{\ast}A\widetilde{U})\|_{F}^{2}-\|\mathcal{U}(\widetilde{U}^{\ast}E\widetilde{U})\|_{F}^{2}\leq\frac{W_{L}(\widetilde{U}^{\ast}A\widetilde{U})}{1+W_{L}(\widetilde{U}^{\ast}A\widetilde{U})}\delta(A)^{2}.
\end{displaymath}
Note that $W_{L}(\widetilde{U}^{\ast}A\widetilde{U})=W_{L}(\widetilde{U}^{\ast}E\widetilde{U})$. An application of (3.2) yields (3.4a).
  
(b) In view of (3.1b), we have
\begin{align*}
\|\mathcal{U}(\widetilde{U}^{\ast}A\widetilde{U})\|_{F}^{2}-\|\mathcal{U}(\widetilde{U}^{\ast}E\widetilde{U})\|_{F}^{2}&=\|\mathcal{U}(\widetilde{U}^{\ast}A\widetilde{U})\|_{F}^{2}-\|\Delta-\mathcal{U}(\widetilde{U}^{\ast}A\widetilde{U})\|_{F}^{2}\\
&\leq\|\mathcal{U}(\widetilde{U}^{\ast}A\widetilde{U})\|_{F}^{2}-\Big(\|\Delta\|_{F}-\|\mathcal{U}(\widetilde{U}^{\ast}A\widetilde{U})\|_{F}\Big)^{2}\\
&=2\|\mathcal{U}(\widetilde{U}^{\ast}A\widetilde{U})\|_{F}\|\Delta\|_{F}-\|\Delta\|_{F}^{2}.
\end{align*}
Using Lemma 3.2 and (3.2), we immediately obtain the estimate (3.4b). \qed

\medskip

We now give an explanation on the quantity $\delta(A)$ involved in (3.4a) and (3.4b). Let $A\in\mathbb{C}^{n\times n}$ be a normal matrix with spectrum $\{\lambda_{i}\}_{i=1}^{n}$. We define
\begin{align*}
\mathbf{1}_{n}:&=(1,\ldots,1),\\
\bm{\lambda}_{R}:&=\big({\rm Re}(\lambda_{1}),\ldots,{\rm Re}(\lambda_{n})\big),\\ \bm{\lambda}_{I}:&=\big({\rm Im}(\lambda_{1}),\ldots,{\rm Im}(\lambda_{n})\big), 
\end{align*}
where ${\rm Re}(\cdot)$ and ${\rm Im}(\cdot)$ denote the real part and the imaginary part of a complex number, respectively. Then we have
\begin{displaymath}
\frac{1}{n}|{\rm tr}(A)|^{2}=\frac{1}{n}\left[(\bm{\lambda}_{R}\cdot\bm{1}_{n})^{2}+(\bm{\lambda}_{I}\cdot\bm{1}_{n})^{2}\right]=|\bm{\lambda}_{R}|^{2}\cos^{2}\theta_{1}+|\bm{\lambda}_{I}|^{2}\cos^{2}\theta_{2},
\end{displaymath}
where 
\begin{displaymath}
\theta_{1}=\arccos\frac{\bm{\lambda}_{R}\cdot\bm{1}_{n}}{|\bm{\lambda}_{R}||\bm{1}_{n}|} \quad \text{and} \quad \theta_{2}=\arccos\frac{\bm{\lambda}_{I}\cdot\bm{1}_{n}}{|\bm{\lambda}_{I}||\bm{1}_{n}|}.
\end{displaymath}
The normality of $A$ implies 
\begin{displaymath}
\|A\|_{F}^{2}=\sum_{i=1}^{n}|\lambda_{i}|^{2}=|\bm{\lambda}_{R}|^{2}+|\bm{\lambda}_{I}|^{2}.
\end{displaymath}
Thus, $\delta(A)$ can be explicitly expressed as
\begin{displaymath}
\delta(A)=\sqrt{|\bm{\lambda}_{R}|^{2}\sin^{2}\theta_{1}+|\bm{\lambda}_{I}|^{2}\sin^{2}\theta_{2}}.
\end{displaymath}

\medskip

A key observation is that the estimates in Lemmas 3.4 and 3.5 are still valid if $W_{L}(\widetilde{U}^{\ast}E\widetilde{U})$ is replaced by $n-1$. Hence, we have the following results.

\medskip

\noindent{\bf Theorem 3.6.} \emph{Under the conditions of Lemma 3.4. Then there is a permutation $\pi$ of $\{1,\ldots,n\}$ such that}
\begin{align*}
\mathbb{D}_{2}&\leq\sqrt{\|E\|_{F}^{2}+(n-1)\delta(E)^{2}},\tag{3.5a}\\
\mathbb{D}_{2}&\leq\sqrt{\|E\|_{F}^{2}+\sqrt{n}\ \delta(E)\|\Delta\|_{F}},\tag{3.5b}\\
\mathbb{D}_{2}&\leq\sqrt{\|E\|_{F}^{2}+2\delta(E)\|\Delta\|_{F}+\|\Delta\|_{F}^{2}},\tag{3.5c}\\
\mathbb{D}_{2}&\leq\sqrt{\|E\|_{F}^{2}+2\sqrt{n-1}\ \delta(E)\|\Delta\|_{F}-\|\Delta\|_{F}^{2}},\tag{3.5d}\\
\mathbb{D}_{2}&\leq\sqrt{\|E\|_{F}^{2}+\frac{n-1}{n}\delta(A)^{2}},\tag{3.5e}\\
\mathbb{D}_{2}&\leq\sqrt{\|E\|_{F}^{2}+2\sqrt{\frac{n-1}{n}}\delta(A)\|\Delta\|_{F}-\|\Delta\|_{F}^{2}},\tag{3.5f}
\end{align*}
\emph{where $\Delta=\mathcal{U}(\widetilde{U}^{\ast}\widetilde{A}\widetilde{U})$ with $\widetilde{U}\in\mathscr{U}_{n}(\widetilde{A})$.}

\medskip

\noindent{\bf Remark 3.7.}  For a given normal matrix $A$, the upper bounds in (3.5a) and (3.5e) only depend on the perturbation $E$. They may not reduce to the Hoffman--Wielandt’s bound $\|E\|_{F}$, even if the perturbed matrix $\widetilde{A}$ is normal. Nevertheless, the upper bounds in (3.5b)--(3.5d) and (3.5f) reduce to $\|E\|_{F}$ when $\widetilde{A}$ is normal. In other words, the estimates (3.5b)--(3.5d) and (3.5f) have extended the Hoffman--Wielandt theorem.

\bigskip

\noindent{\bf Remark 3.8.}  Although $\|\Delta\|_{F}$ can be explicitly expressed as
\begin{displaymath} \|\Delta\|_{F}=\bigg(\|\widetilde{A}\|_{F}^{2}-\sum_{i=1}^{n}\big|\widetilde{\lambda}_{i}\big|^{2}\bigg)^{1\over 2},
\end{displaymath}
the quantity $\|\Delta\|_{F}$ is uncomputable in general because the spectrum of $\widetilde{A}$ is unknown. Hence, a natural question is how to effectively estimate  $\|\Delta\|_{F}$. Here we mention an applicable upper bound for $\|\Delta\|_{F}$ derived by Henrici~[8], that is,
\begin{align*}
\|\Delta\|_{F}\leq\left(\frac{n^{3}-n}{12}\right)^{1\over 4}\sqrt{\|\widetilde{A}\widetilde{A}^{\ast}-\widetilde{A}^{\ast}\widetilde{A}\|_{F}}.\tag{3.6}
\end{align*}
There is also a lower bound for $\|\Delta\|_{F}$ established by Sun [9], that is, 
\begin{align*}
\|\Delta\|_{F}\geq\left(\|\widetilde{A}\|_{F}^{2}-\sqrt{\|\widetilde{A}\|_{F}^{4}-\frac{1}{2}\|\widetilde{A}\widetilde{A}^{\ast}-\widetilde{A}^{\ast}\widetilde{A}\|_{F}^{2}}\right)^{1\over 2}.\tag{3.7}
\end{align*}
The inequalities (3.6) and (3.7) justify that $\|\Delta\|_{F}$ can be viewed as a measure of the non-normality of $\widetilde{A}$. From (3.6) and (3.7), we see that $\widetilde{A}$ is normal if and only if $\|\Delta\|_{F}=0$. 

\medskip

\noindent{\bf Remark 3.9.} Define \begin{displaymath}
\mathbb{D}_{\infty}:=\max_{1\leq i\leq n}\big|\widetilde{\lambda}_{\pi(i)}-\lambda_{i}\big|.
\end{displaymath}
Using $\mathbb{D}_{\infty}\leq\mathbb{D}_{2}$ and Theorem 3.6, we can readily get the corresponding estimates for $\mathbb{D}_{\infty}$.

\bigskip

{\bf 3.3. The estimates based on block decomposition.} Assume that $\widetilde{U}_{1}\in\mathscr{U}_{n}$ satisfies
\begin{align*}
\widetilde{A}=\widetilde{U}_{1}{\rm diag}\big(\widetilde{A}_{1},\ldots,\widetilde{A}_{s}\big)\widetilde{U}_{1}^{\ast}\tag{3.8}
\end{align*}
for some positive integer $s$, where each $\widetilde{A}_{i}\in\mathbb{C}^{n_{i}\times n_{i}}$ is upper triangular and $\sum_{i=1}^{s}n_{i}=n$. In particular, if $s=1$, (3.8) is the Schur’s decomposition of $\widetilde{A}$; if $s=n$, (3.8) implies that $\widetilde{A}$ is normal. We then have
\begin{displaymath}
\widetilde{U}_{1}^{\ast}\widetilde{A}\widetilde{U}_{1}={\rm diag}\big(\widetilde{A}_{1},\ldots,\widetilde{A}_{s}\big)=\widetilde{\Lambda}+\Delta_{1},
\end{displaymath}
where $\widetilde{\Lambda}={\rm diag}\big(\mathcal{D}(\widetilde{A}_{1}),\ldots,\mathcal{D}(\widetilde{A}_{s})\big)$ and $\Delta_{1}={\rm diag}\big(\mathcal{U}(\widetilde{A}_{1}),\ldots,\mathcal{U}(\widetilde{A}_{s})\big)$.

Let $\widetilde{U}_{1}^{\ast}A\widetilde{U}_{1}$ and $\widetilde{U}_{1}^{\ast}E\widetilde{U}_{1}$ be partitioned as the block forms that coincide with $\widetilde{U}_{1}^{\ast}\widetilde{A}\widetilde{U}_{1}$. Set $\widetilde{U}_{1}^{\ast}A\widetilde{U}_{1}=(\widehat{A}_{ij})_{s\times s}$ and $\widetilde{U}_{1}^{\ast}E\widetilde{U}_{1}=(\widehat{E}_{ij})_{s\times s}$. Due to $\widetilde{U}_{1}^{\ast}A\widetilde{U}_{1}+\widetilde{U}_{1}^{\ast}E\widetilde{U}_{1}={\rm diag}\big(\widetilde{A}_{1},\ldots,\widetilde{A}_{s}\big)$, it follows that
\begin{align*}
\widehat{A}_{ii}+\widehat{E}_{ii}=\widetilde{A}_{i}, \quad \forall i=1,\ldots,s.\tag{3.9}
\end{align*}
By the Hoffman--Wielandt theorem, we have that there exists a permutation $\pi$ of $\{1,\ldots,n\}$ such that
\begin{displaymath}
\mathbb{D}_{2}\leq\|\widetilde{\Lambda}-\widetilde{U}_{1}^{\ast}A\widetilde{U}_{1}\|_{F}=\|\widetilde{U}_{1}^{\ast}E\widetilde{U}_{1}-\Delta_{1}\|_{F}.
\end{displaymath}
Hence, we have
\begin{align*}
\mathbb{D}_{2}\leq\sqrt{\|\mathcal{D}(\widetilde{U}_{1}^{\ast}E\widetilde{U}_{1})\|_{F}^{2}+\|\mathcal{L}(\widetilde{U}_{1}^{\ast}E\widetilde{U}_{1})\|_{F}^{2}+\|\mathcal{U}(\widetilde{U}_{1}^{\ast}E\widetilde{U}_{1})-\Delta_{1}\|_{F}^{2}}.\tag{3.10}
\end{align*}

The following theorem presents another three estimates for $\mathbb{D}_{2}$, which involve the quantity $s$.

\medskip

\noindent{\bf Theorem 3.10.} \emph{Let $A\in\mathbb{C}^{n\times n}$ be normal, and let $\widetilde{A}=A+E$ with decomposition {\rm (3.8)}, where $E\in\mathbb{C}^{n\times n}$ is an arbitrary perturbation. Assume that  the spectra of $A$ and $\widetilde{A}$ are $\{\lambda_{i}\}_{i=1}^{n}$ and $\{\widetilde{\lambda}_{i}\}_{i=1}^{n}$, respectively. Then there exists a permutation $\pi$ of $\{1,\ldots,n\}$ such that}
\begin{align*}
\mathbb{D}_{2}&\leq\sqrt{\|E\|_{F}^{2}+(n-s)\delta(E)^{2}},\tag{3.11a}\\
\mathbb{D}_{2}&\leq\sqrt{\|E\|_{F}^{2}+\sqrt{n-s+1}\ \delta(E)\|\Delta_{1}\|_{F}},\tag{3.11b}\\
\mathbb{D}_{2}&\leq\sqrt{\|E\|_{F}^{2}+2\sqrt{n-s}\ \delta(E)\|\Delta_{1}\|_{F}-\|\Delta_{1}\|_{F}^{2}},\tag{3.11c}
\end{align*}
\emph{where $\Delta_{1}=\mathcal{U}(\widetilde{U}_{1}^{\ast}\widetilde{A}\widetilde{U}_{1})$.}

\medskip

\noindent{\bf Proof.} (a) By (3.9), we have 
\begin{align*}
\|\mathcal{U}(\widetilde{U}_{1}^{\ast}E\widetilde{U}_{1})-\Delta_{1}\|_{F}^{2}&=\sum_{1\leq i<j\leq s}\|\widehat{E}_{ij}\|_{F}^{2}+\sum_{i=1}^{s}\|\mathcal{U}(\widehat{E}_{ii})-\mathcal{U}(\widetilde{A}_{i})\|_{F}^{2}\tag{3.12a}\\
&\leq\|\mathcal{U}(\widetilde{U}_{1}^{\ast}E\widetilde{U}_{1})\|_{F}^{2}+\sum_{i=1}^{s}\|\mathcal{U}(\widehat{A}_{ii})\|_{F}^{2}.\tag{3.12b}
\end{align*}
Let $(\mathcal{U}(\widehat{A}_{ii}))_{(j)}$ be the $j$-th row of $\mathcal{U}(\widehat{A}_{ii})$, and let 
\begin{displaymath}
\|(\mathcal{U}(\widehat{A}_{pp}))_{(q)}\|_{F}=\max\limits_{i,j}\|(\mathcal{U}(\widehat{A}_{ii}))_{(j)}\|_{F}
\end{displaymath}
for some $1\leq p\leq s$ and $1\leq q\leq n_{p}$. Then we have
\begin{displaymath}
\sum_{i=1}^{s}\|\mathcal{U}(\widehat{A}_{ii})\|_{F}^{2}\leq\sum_{i=1}^{s}(n_{i}-1)\|(\mathcal{U}(\widehat{A}_{pp}))_{(q)}\|_{F}^{2}\leq(n-s)\|(\mathcal{U}(\widetilde{U}_{1}^{\ast}A\widetilde{U}_{1}))_{(k)}\|_{F}^{2},
\end{displaymath}
where $k=\sum_{i=1}^{p-1}n_{i}+q$ (if $p=1$, we set $k=q$). Using Lemma 2.4, we have
\begin{align*}
\sum_{i=1}^{s}\|\mathcal{U}(\widehat{A}_{ii})\|_{F}^{2}\leq(n-s)\|\mathcal{L}(\widetilde{U}_{1}^{\ast}A\widetilde{U}_{1})\|_{F}^{2}=(n-s)\|\mathcal{L}(\widetilde{U}_{1}^{\ast}E\widetilde{U}_{1})\|_{F}^{2}.\tag{3.13}
\end{align*}
Combining (3.12b) and (3.13), we obtain
\begin{displaymath}
\|\mathcal{U}(\widetilde{U}_{1}^{\ast}E\widetilde{U}_{1})-\Delta_{1}\|_{F}^{2}\leq\|\mathcal{U}(\widetilde{U}_{1}^{\ast}E\widetilde{U}_{1})\|_{F}^{2}+(n-s)\|\mathcal{L}(\widetilde{U}_{1}^{\ast}E\widetilde{U}_{1})\|_{F}^{2}.
\end{displaymath}
Due to (3.10) and $\|\mathcal{L}(\widetilde{U}_{1}^{\ast}E\widetilde{U}_{1})\|_{F}\leq\delta(E)$, it follows that (3.11a) holds.

(b) From (3.9) and (3.12a), we have
\begin{align*}
\|\mathcal{U}(\widetilde{U}_{1}^{\ast}E\widetilde{U}_{1})-\Delta_{1}\|_{F}^{2}=\|\mathcal{U}(\widetilde{U}_{1}^{\ast}E\widetilde{U}_{1})\|_{F}^{2}+\sum_{i=1}^{s}\left(\|\mathcal{U}(\widehat{A}_{ii})\|_{F}^{2}-\|\mathcal{U}(\widehat{E}_{ii})\|_{F}^{2}\right).\tag{3.14}
\end{align*}
In view of (3.9) and the triangle inequality, we immediately obtain
\begin{displaymath}
\|\mathcal{U}(\widehat{A}_{ii})\|_{F}-\|\mathcal{U}(\widehat{E}_{ii})\|_{F}\leq\|\mathcal{U}(\widetilde{A}_{i})\|_{F}, \quad \forall i=1,\ldots,s.
\end{displaymath}
Hence,
\begin{displaymath}
\|\mathcal{U}(\widetilde{U}_{1}^{\ast}E\widetilde{U}_{1})-\Delta_{1}\|_{F}^{2}\leq\|\mathcal{U}(\widetilde{U}_{1}^{\ast}E\widetilde{U}_{1})\|_{F}^{2}+\sum_{i=1}^{s}\|\mathcal{U}(\widetilde{A}_{i})\|_{F}\Big(\|\mathcal{U}(\widehat{A}_{ii})\|_{F}+\|\mathcal{U}(\widehat{E}_{ii})\|_{F}\Big).
\end{displaymath}
By the Cauchy--Schwarz’s inequality, we have
\begin{align*}
\sum_{i=1}^{s}\|\mathcal{U}(\widetilde{A}_{i})\|_{F}\|\mathcal{U}(\widehat{A}_{ii})\|_{F}&\leq\|\Delta_{1}\|_{F}\bigg(\sum_{i=1}^{s}\|\mathcal{U}(\widehat{A}_{ii})\|_{F}^{2}\bigg)^{1\over 2},\\
\sum_{i=1}^{s}\|\mathcal{U}(\widetilde{A}_{i})\|_{F}\|\mathcal{U}(\widehat{E}_{ii})\|_{F}&\leq\|\Delta_{1}\|_{F}\bigg(\sum_{i=1}^{s}\|\mathcal{U}(\widehat{E}_{ii})\|_{F}^{2}\bigg)^{1\over 2}.
\end{align*}
Since
\begin{displaymath}
\sum_{i=1}^{s}\|\mathcal{U}(\widehat{A}_{ii})\|_{F}^{2}\leq(n-s)\|\mathcal{L}(\widetilde{U}_{1}^{\ast}E\widetilde{U}_{1})\|_{F}^{2} \quad \text{and} \quad \sum_{i=1}^{s}\|\mathcal{U}(\widehat{E}_{ii})\|_{F}^{2}\leq\|\mathcal{U}(\widetilde{U}_{1}^{\ast}E\widetilde{U}_{1})\|_{F}^{2},
\end{displaymath}
it follows that
\begin{displaymath}
\|\mathcal{U}(\widetilde{U}_{1}^{\ast}E\widetilde{U}_{1})-\Delta_{1}\|_{F}^{2}\leq\|\mathcal{U}(\widetilde{U}_{1}^{\ast}E\widetilde{U}_{1})\|_{F}^{2}+\Big(\sqrt{n-s}\|\mathcal{L}(\widetilde{U}_{1}^{\ast}E\widetilde{U}_{1})\|_{F}+\|\mathcal{U}(\widetilde{U}_{1}^{\ast}E\widetilde{U}_{1})\|_{F}\Big)\|\Delta_{1}\|_{F}.
\end{displaymath}
Using the Cauchy--Schwarz’s inequality and Lemma 3.1, we obtain
\begin{displaymath} 
\sqrt{n-s}\|\mathcal{L}(\widetilde{U}_{1}^{\ast}E\widetilde{U}_{1})\|_{F}+\|\mathcal{U}(\widetilde{U}_{1}^{\ast}E\widetilde{U}_{1})\|_{F}\leq\sqrt{n-s+1}\ \delta(E).
\end{displaymath}
Thus, 
\begin{displaymath}
\|\mathcal{U}(\widetilde{U}_{1}^{\ast}E\widetilde{U}_{1})-\Delta_{1}\|_{F}^{2}\leq\|\mathcal{U}(\widetilde{U}_{1}^{\ast}E\widetilde{U}_{1})\|_{F}^{2}+\sqrt{n-s+1}\ \delta(E)\|\Delta_{1}\|_{F}.
\end{displaymath}
Then, the estimate (3.11b) follows immediately from (3.10).

(c) By (3.9) and the triangle inequality, we have
\begin{align*}
\sum_{i=1}^{s}\Big(\|\mathcal{U}(\widehat{A}_{ii})\|_{F}^{2}-\|\mathcal{U}(\widehat{E}_{ii})\|_{F}^{2}\Big)&=\sum_{i=1}^{s}\Big(\|\mathcal{U}(\widehat{A}_{ii})\|_{F}^{2}-\|\mathcal{U}(\widetilde{A}_{i})-\mathcal{U}(\widehat{A}_{ii})\|_{F}^{2}\Big)\\
&\leq\sum_{i=1}^{s}\bigg[\|\mathcal{U}(\widehat{A}_{ii})\|_{F}^{2}-\Big(\|\mathcal{U}(\widetilde{A}_{i})\|_{F}-\|\mathcal{U}(\widehat{A}_{ii})\|_{F}\Big)^{2}\bigg]\\
&=2\sum_{i=1}^{s}\|\mathcal{U}(\widetilde{A}_{i})\|_{F}\|\mathcal{U}(\widehat{A}_{ii})\|_{F}-\|\Delta_{1}\|_{F}^{2}.
\end{align*}
Using the Cauchy--Schwarz’s inequality, we obtain
\begin{align*}
\sum_{i=1}^{s}\left(\|\mathcal{U}(\widehat{A}_{ii})\|_{F}^{2}-\|\mathcal{U}(\widehat{E}_{ii})\|_{F}^{2}\right)\leq 2\bigg(\sum_{i=1}^{s}\|\mathcal{U}(\widehat{A}_{ii})\|_{F}^{2}\bigg)^{1\over 2}\|\Delta_{1}\|_{F}-\|\Delta_{1}\|_{F}^{2}.\tag{3.15}
\end{align*}
In view of (3.13)--(3.15), we arrive at
\begin{displaymath}
\|\mathcal{U}(\widetilde{U}_{1}^{\ast}E\widetilde{U}_{1})-\Delta_{1}\|_{F}^{2}\leq\|\mathcal{U}(\widetilde{U}_{1}^{\ast}E\widetilde{U}_{1})\|_{F}^{2}+2\sqrt{n-s}\|\mathcal{L}(\widetilde{U}_{1}^{\ast}E\widetilde{U}_{1})\|_{F}\|\Delta_{1}\|_{F}-\|\Delta_{1}\|_{F}^{2}.
\end{displaymath}
From (3.10) and $\|\mathcal{L}(\widetilde{U}_{1}^{\ast}E\widetilde{U}_{1})\|_{F}\leq\delta(E)$, we deduce that (3.11c) holds. \qed

\medskip

\noindent{\bf Remark 3.11.} For any $\widetilde{A}\in\mathbb{C}^{n\times n}$, the decomposition (3.8) is always valid for $s=1$. If $s=1$, (3.11a)--(3.11c) reduce to (3.5a), (3.5b), and (3.5d), respectively; otherwise, the upper bounds in (3.11a)--(3.11c) are smaller.

\medskip

\noindent{\bf Remark 3.12.} We remark that the upper bounds for $\mathbb{D}_{2}$ derived in this section are all absolute type perturbation bounds. As discussed in [6], we can apply our estimates for $\|\widetilde{U}^{\ast}E\widetilde{U}-\Delta\|_{F}$ and $\|\widetilde{U}_{1}^{\ast}E\widetilde{U}_{1}-\Delta_{1}\|_{F}$ to derive corresponding relative type perturbation bounds. For more theories about relative perturbation bounds of spectrum, we refer to [12--15] and the references therein. 

\medskip
\bigskip

\noindent{\bf \large 4. Perturbation bounds for the spectrum of a Hermitian matrix}

\medskip

Clearly, the estimates established in Section 3 are applicable for Hermitian matrices as well. Nevertheless, Hermitian matrices possess some special properties, which can yield more accurate estimates. In this section, we present several new perturbation bounds for the spectrum of a Hermitian matrix. 

Let $A\in\mathbb{C}^{n\times n}$ be Hermitian, and let its perturbed matrix $\widetilde{A}\in\mathbb{C}^{n\times n}$ be decomposed as
\begin{align*}
\widetilde{A}=\widetilde{U}\big(\widetilde{\Lambda}+\Delta\big)\widetilde{U}^{\ast},\tag{4.1}
\end{align*}
where $\widetilde{U}\in\mathscr{U}_{n}$, $\widetilde{\Lambda}={\rm diag}\big(\widetilde{\lambda}_{1},\ldots,\widetilde{\lambda}_{n}\big)$, and $\Delta$ is strictly upper triangular.

Using the same argument as in subsection 3.2, we can obtain
\begin{align*}
\|\widetilde{U}^{\ast}E\widetilde{U}-\Delta\|_{F}^{2}=\|E\|_{F}^{2}+\|\mathcal{U}(\widetilde{U}^{\ast}A\widetilde{U})\|_{F}^{2}-\|\mathcal{U}(\widetilde{U}^{\ast}E\widetilde{U})\|_{F}^{2}.\tag{4.2}
\end{align*}
Since $\widetilde{U}^{\ast}A\widetilde{U}$ is Hermitian, by (3.1a), we have
\begin{displaymath}
\|\mathcal{U}(\widetilde{U}^{\ast}A\widetilde{U})\|_{F}^{2}=\|\mathcal{L}(\widetilde{U}^{\ast}A\widetilde{U})\|_{F}^{2}=\|\mathcal{L}(\widetilde{U}^{\ast}E\widetilde{U})\|_{F}^{2}.
\end{displaymath}
Hence,
\begin{align*}
\|\widetilde{U}^{\ast}E\widetilde{U}-\Delta\|_{F}^{2}=\|E\|_{F}^{2}+\|\mathcal{L}(\widetilde{U}^{\ast}E\widetilde{U})\|_{F}^{2}-\|\mathcal{U}(\widetilde{U}^{\ast}E\widetilde{U})\|_{F}^{2}.\tag{4.3}
\end{align*}

In order to derive perturbation bounds for the spectrum of a Hermitian matrix, we need the following estimates for $\|\widetilde{U}^{\ast}E\widetilde{U}-\Delta\|_{F}$.

\medskip

\noindent{\bf Lemma 4.1.} \emph{Let $A\in\mathbb{C}^{n\times n}$ be Hermitian, and let $\widetilde{A}=A+E$, where $E\in\mathbb{C}^{n\times n}$ is an arbitrary perturbation. Then}
\begin{align*}
\|\widetilde{U}^{\ast}E\widetilde{U}-\Delta\|_{F}&\leq\sqrt{\|E\|_{F}^{2}+\delta(E)^{2}},\tag{4.4a}\\
\|\widetilde{U}^{\ast}E\widetilde{U}-\Delta\|_{F}&\leq\sqrt{\|E\|_{F}^{2}+\sqrt{2}\delta(E)\|\Delta\|_{F}},\tag{4.4b}\\
\|\widetilde{U}^{\ast}E\widetilde{U}-\Delta\|_{F}&\leq\sqrt{\|E\|_{F}^{2}+2\delta(E)\|\Delta\|_{F}-\|\Delta\|_{F}^{2}},\tag{4.4c}\\
\|\widetilde{U}^{\ast}E\widetilde{U}-\Delta\|_{F}&\leq\sqrt{\|E\|_{F}^{2}+\frac{1}{2}\delta(A)^{2}},\tag{4.4d}\\
\|\widetilde{U}^{\ast}E\widetilde{U}-\Delta\|_{F}&\leq\sqrt{\|E\|_{F}^{2}+\sqrt{2}\delta(A)\|\Delta\|_{F}-\|\Delta\|_{F}^{2}},\tag{4.4e}
\end{align*}
\emph{where $\Delta=\mathcal{U}(\widetilde{U}^{\ast}\widetilde{A}\widetilde{U})$ with $\widetilde{U}\in\mathscr{U}_{n}(\widetilde{A})$.}

\medskip

\noindent{\bf Proof.} 
(a) Using (4.3) and Lemma 3.1, we obtain
\begin{align*}
\|\widetilde{U}^{\ast}E\widetilde{U}-\Delta\|_{F}&\leq \sqrt{\|E\|_{F}^{2}+\|\mathcal{L}(\widetilde{U}^{\ast}E\widetilde{U})\|_{F}^{2}+\|\mathcal{U}(\widetilde{U}^{\ast}E\widetilde{U})\|_{F}^{2}}\\
&\leq \sqrt{\|E\|_{F}^{2}+\delta(E)^{2}}.
\end{align*}

(b) From (3.1b), we have that $\|\mathcal{U}(\widetilde{U}^{\ast}A\widetilde{U})\|_{F}-\|\mathcal{U}(\widetilde{U}^{\ast}E\widetilde{U})\|_{F}\leq\|\Delta\|_{F}$. In addition,
\begin{align*}
\|\mathcal{U}(\widetilde{U}^{\ast}A\widetilde{U})\|_{F}+\|\mathcal{U}(\widetilde{U}^{\ast}E\widetilde{U})\|_{F}&=\|\mathcal{L}(\widetilde{U}^{\ast}E\widetilde{U})\|_{F}+\|\mathcal{U}(\widetilde{U}^{\ast}E\widetilde{U})\|_{F}\\
&\leq\sqrt{2}\sqrt{\|\mathcal{L}(\widetilde{U}^{\ast}E\widetilde{U})\|_{F}^{2}+\|\mathcal{U}(\widetilde{U}^{\ast}E\widetilde{U})\|_{F}^{2}}\\
&\leq\sqrt{2}\delta(E).
\end{align*}
An application of (4.2) yields the inequality (4.4b).

(c) By (4.2), (3.1b), and the triangle inequality, we have
\begin{align*}
\|\widetilde{U}^{\ast}E\widetilde{U}-\Delta\|_{F}^{2}&\leq\|E\|_{F}^{2}+\|\mathcal{U}(\widetilde{U}^{\ast}A\widetilde{U})\|_{F}^{2}-\Big(\|\Delta\|_{F}-\|\mathcal{U}(\widetilde{U}^{\ast}A\widetilde{U})\|_{F}\Big)^{2}\\
&=\|E\|_{F}^{2}+2\|\mathcal{U}(\widetilde{U}^{\ast}A\widetilde{U})\|_{F}\|\Delta\|_{F}-\|\Delta\|_{F}^{2}\\
&=\|E\|_{F}^{2}+2\|\mathcal{L}(\widetilde{U}^{\ast}E\widetilde{U})\|_{F}\|\Delta\|_{F}-\|\Delta\|_{F}^{2}\\
&\leq\|E\|_{F}^{2}+2\sqrt{\|\mathcal{L}(\widetilde{U}^{\ast}E\widetilde{U})\|_{F}^{2}+\|\mathcal{U}(\widetilde{U}^{\ast}E\widetilde{U})\|_{F}^{2}}\|\Delta\|_{F}-\|\Delta\|_{F}^{2}\\
&\leq\|E\|_{F}^{2}+2\delta(E)\|\Delta\|_{F}-\|\Delta\|_{F}^{2},
\end{align*}
which gives the inequality (4.4c).

(d) Because $\widetilde{U}^{\ast}A\widetilde{U}$ is Hermitian, by (4.2) and Lemma 3.1, we have
\begin{align*}
\|\widetilde{U}^{\ast}E\widetilde{U}-\Delta\|_{F}^{2}&=\|E\|_{F}^{2}+\|\mathcal{U}(\widetilde{U}^{\ast}A\widetilde{U})\|_{F}^{2}-\|\mathcal{U}(\widetilde{U}^{\ast}E\widetilde{U})\|_{F}^{2}\\
&\leq\|E\|_{F}^{2}+\frac{1}{2}\big(\|\mathcal{L}(\widetilde{U}^{\ast}A\widetilde{U})\|_{F}^{2}+\|\mathcal{U}(\widetilde{U}^{\ast}A\widetilde{U})\|_{F}^{2}\big)\\
&\leq\|E\|_{F}^{2}+\frac{1}{2}\delta(A)^{2},
\end{align*}
which means the inequality (4.4d).

(e) Based on the derivation in (c), we have
\begin{align*}
\|\widetilde{U}^{\ast}E\widetilde{U}-\Delta\|_{F}^{2}&\leq\|E\|_{F}^{2}+2\|\mathcal{U}(\widetilde{U}^{\ast}A\widetilde{U})\|_{F}\|\Delta\|_{F}-\|\Delta\|_{F}^{2}\\
&=\|E\|_{F}^{2}+\big(\|\mathcal{U}(\widetilde{U}^{\ast}A\widetilde{U})\|_{F}+\|\mathcal{L}(\widetilde{U}^{\ast}A\widetilde{U})\|_{F}\big)\|\Delta\|_{F}-\|\Delta\|_{F}^{2}\\
&\leq\|E\|_{F}^{2}+\sqrt{2}\sqrt{\|\mathcal{U}(\widetilde{U}^{\ast}A\widetilde{U})\|_{F}^{2}+\|\mathcal{L}(\widetilde{U}^{\ast}A\widetilde{U})\|_{F}^{2}}\|\Delta\|_{F}-\|\Delta\|_{F}^{2}\\
&\leq\|E\|_{F}^{2}+\sqrt{2}\delta(A)\|\Delta\|_{F}-\|\Delta\|_{F}^{2},
\end{align*}
which leads to the inequality (4.4e). \qed

\medskip

Since both $\widetilde{\Lambda}$ and $\widetilde{U}^{\ast}A\widetilde{U}$ are normal, by the Hoffman--Wielandt theorem, we get that there exists a permutation $\pi$ of $\{1,\ldots,n\}$ such that
\begin{align*}
\mathbb{D}_{2}\leq\|\widetilde{\Lambda}-\widetilde{U}^{\ast}A\widetilde{U}\|_{F}=\|\widetilde{U}^{\ast}E\widetilde{U}-\Delta\|_{F}.\tag{4.5}
\end{align*}
An application of Lemma 4.1 yields the following theorem.

\medskip

\noindent{\bf Theorem 4.2.} \emph{Let $A\in\mathbb{C}^{n\times n}$ be a Hermitian matrix with spectrum $\{\lambda_{i}\}_{i=1}^{n}$, and let $\widetilde{A}=A+E$ with spectrum $\{\widetilde{\lambda}_{i}\}_{i=1}^{n}$, where $E\in\mathbb{C}^{n\times n}$ is an arbitrary perturbation. Then there is a permutation $\pi$ of $\{1,\ldots,n\}$ such that}
\begin{align*}
\mathbb{D}_{2}&\leq\sqrt{\|E\|_{F}^{2}+\delta(E)^{2}},\tag{4.6a}\\
\mathbb{D}_{2}&\leq\sqrt{\|E\|_{F}^{2}+\sqrt{2}\delta(E)\|\Delta\|_{F}},\tag{4.6b}\\
\mathbb{D}_{2}&\leq\sqrt{\|E\|_{F}^{2}+2\delta(E)\|\Delta\|_{F}-\|\Delta\|_{F}^{2}},\tag{4.6c}\\
\mathbb{D}_{2}&\leq\sqrt{\|E\|_{F}^{2}+\frac{1}{2}\delta(A)^{2}},\tag{4.6d}\\
\mathbb{D}_{2}&\leq\sqrt{\|E\|_{F}^{2}+\sqrt{2}\delta(A)\|\Delta\|_{F}-\|\Delta\|_{F}^{2}},\tag{4.6e}
\end{align*}
\emph{where $\Delta=\mathcal{U}(\widetilde{U}^{\ast}\widetilde{A}\widetilde{U})$ with $\widetilde{U}\in\mathscr{U}_{n}(\widetilde{A})$.}

\medskip

We next give two useful estimates for $\|\Delta\|_{F}$, which are sharper than the existing results in [10, Lemma 2.1 and Corollary 2.5].

\medskip

\noindent{\bf Theorem 4.3.} \emph{Let $A\in\mathbb{C}^{n\times n}$ be Hermitian, and let $\widetilde{A}=A+E \ (\widetilde{A}\neq 0)$, where $E\in\mathbb{C}^{n\times n}$ is a perturbation. Assume that the Schur's decomposition {\rm (4.1)} satisfies that $\big|\widetilde{\lambda}_{1}\big|\geq\cdots\geq\big|\widetilde{\lambda}_{n}\big|$. Then}
\begin{align*}
\|\Delta\|_{F}&\leq\frac{1}{\sqrt{2}}\bigg(\|\widetilde{A}-\widetilde{A}^{\ast}\|_{F}^{2}-\frac{1}{{\rm rank}(\widetilde{A})}\big|{\rm tr}(\widetilde{A}-\widetilde{A}^{\ast})\big|^{2}\bigg)^{1\over 2}\leq\frac{1}{\sqrt{2}}\delta\big(\widetilde{A}-\widetilde{A}^{\ast}\big),\\
\|\Delta\|_{F}&\leq\frac{1}{\sqrt{2}}\bigg(\|E-E^{\ast}\|_{F}^{2}-\frac{1}{{\rm rank}(\widetilde{A})}\big|{\rm tr}(E-E^{\ast})\big|^{2}\bigg)^{1\over 2}\leq\frac{1}{\sqrt{2}}\delta(E-E^{\ast}).
\end{align*}

\noindent{\bf Proof.} By decomposition (4.1), we have 
\begin{displaymath}
\widetilde{A}-\widetilde{A}^{\ast}=\widetilde{U}(\widetilde{\Lambda}+\Delta-\widetilde{\Lambda}^{\ast}-\Delta^{\ast})\widetilde{U}^{\ast},
\end{displaymath}
which gives 
\begin{displaymath}
\|\widetilde{A}-\widetilde{A}^{\ast}\|_{F}^{2}=\|\widetilde{\Lambda}-\widetilde{\Lambda}^{\ast}\|_{F}^{2}+\|\Delta-\Delta^{\ast}\|_{F}^{2}=\sum_{i=1}^{\rm rank{(\widetilde{A})}}\Big|\widetilde{\lambda}_{i}-\overline{\widetilde{\lambda}}_{i}\Big|^{2}+2\|\Delta\|_{F}^{2}.
\end{displaymath}
Using the Cauchy--Schwarz’s inequality, we obtain
\begin{align*}
\|\widetilde{A}-\widetilde{A}^{\ast}\|_{F}^{2}&\geq\frac{1}{{\rm rank}(\widetilde{A})}\Bigg(\sum_{i=1}^{\rm rank{(\widetilde{A})}}\Big|\widetilde{\lambda}_{i}-\overline{\widetilde{\lambda}}_{i}\Big|\Bigg)^{2}+2\|\Delta\|_{F}^{2}\\
&\geq\frac{1}{{\rm rank}(\widetilde{A})}\big|{\rm tr}(\widetilde{A}-\widetilde{A}^{\ast})\big|^{2}+2\|\Delta\|_{F}^{2}.
\end{align*}
Consequently, we arrive at
\begin{displaymath}
\|\Delta\|_{F}\leq\frac{1}{\sqrt{2}}\bigg(\|\widetilde{A}-\widetilde{A}^{\ast}\|_{F}^{2}-\frac{1}{{\rm rank}(\widetilde{A})}\big|{\rm tr}(\widetilde{A}-\widetilde{A}^{\ast})\big|^{2}\bigg)^{1\over 2}.
\end{displaymath}
The second inequality follows immediately from  $\widetilde{A}-\widetilde{A}^{\ast}=E-E^{\ast}$. \qed

\medskip

\noindent{\bf Example 4.4.} In order to illustrate the estimates in Theorem 4.2, we
give an example as follows:
\begin{displaymath}
A=\begin{pmatrix}
0 & 1 & 0 & \cdots & 0\\
1 & 0 & 0 & \cdots & 0\\
0 & 0 & 1 & \cdots & 0\\
\vdots & \vdots & \vdots & \ddots & \vdots\\
0 & 0 & 0 & \cdots & 1\\
\end{pmatrix}\in\mathbb{R}^{n\times n} \quad \text{and} \quad E=\begin{pmatrix}
0 & 0 & 0 & \cdots & 0\\
-1 & 0 & 0 & \cdots & 0\\
0 & 0 & -1 & \cdots & 0\\
\vdots & \vdots & \vdots & \ddots & \vdots\\
0 & 0 & 0 & \cdots & -1\\
\end{pmatrix}\in\mathbb{R}^{n\times n} \ (n\geq 3).
\end{displaymath}
We can easily get that
$\mathbb{D}_{2}=\sqrt{n}$ for any permutation $\pi$ of $\{1,\ldots,n\}$. Direct calculations yield 
\begin{displaymath}
\|A\|_{F}^{2}=n, \quad |{\rm tr}(A)|^{2}=(n-2)^{2},\quad \|E\|_{F}^{2}=n-1, \quad |{\rm tr}(E)|^{2}=(n-2)^{2}, \quad \|\Delta\|_{F}=1.
\end{displaymath}
We then have
\begin{displaymath} \delta(E)=\sqrt{3-\frac{4}{n}}, \quad \delta(A)=2\sqrt{1-\frac{1}{n}}.
\end{displaymath}
Hence, the upper bounds in (4.6a)--(4.6e) are
\begin{displaymath}
\sqrt{n-\frac{4}{n}+2}, \ \
\sqrt{n+\sqrt{6-\frac{8}{n}}-1},\ \  \sqrt{n+2\sqrt{3-\frac{4}{n}}-2},\ \ \sqrt{n-\frac{2}{n}+1}, \ \ \sqrt{n+2\sqrt{2-\frac{2}{n}}-2},
\end{displaymath}
respectively.

\medskip

\noindent{\bf Remark 4.5.} From (4.5), we can readily obtain an upper bound $\|E\|_{F}+\|\Delta\|_{F}$ via the triangle inequality. Due to the fact that $\delta(E)\leq\|E\|_{F}$, it follows that
\begin{align*}
\sqrt{\|E\|_{F}^{2}+\sqrt{2}\delta(E)\|\Delta\|_{F}}&\leq\|E\|_{F}+\|\Delta\|_{F},\\
\sqrt{\|E\|_{F}^{2}+2\delta(E)\|\Delta\|_{F}-\|\Delta\|_{F}^{2}}&\leq\|E\|_{F}+\|\Delta\|_{F}.
\end{align*}
Hence, the upper bounds in (4.6b) and (4.6c) are non-trivial. It is worth mentioning that the upper bound in (4.6c) also satisfies
\begin{displaymath}
\sqrt{\|E\|_{F}^{2}+2\delta(E)\|\Delta\|_{F}-\|\Delta\|_{F}^{2}}\leq\sqrt{2}\|E\|_{F},
\end{displaymath}
which reveals that (4.6c) is sharper than (1.6).

\medskip

\noindent{\bf Remark 4.6.} For a given Hermitian matrix $A\in\mathbb{C}^{n\times n}$, the upper bounds in (4.6a) and (4.6d) only depend on the perturbation $E$, while other upper bounds in Theorem 4.2 are related to both $E$ and $\|\Delta\|_{F}$ (Remark 3.8 and Theorem 4.3 have provided several applicable estimates for $\|\Delta\|_{F}$). Clearly, if $\widetilde{A}$ is normal (i.e., $\|\Delta\|_{F}=0$), the upper bounds in (4.6b), (4.6c), and (4.6e) all reduce to the Hoffman--Wielandt’s bound $\|E\|_{F}$. 

\medskip
\bigskip

\noindent{\bf \large 5. Conclusions}

\medskip

In this paper, we have established novel perturbation bounds for the spectrum of a normal matrix (including the case of Hermitian matrices). Some of our estimates improve the existing results in [2, 6, 7, 9, 10]. Moreover, if the perturbed matrix is still normal, the upper bounds involving the ``departure from normality'' of the perturbed matrix reduce to the Hoffman--Wielandt’s bound. Therefore, these estimates have generalized the classical Hoffman--Wielandt theorem.

\medskip
\bigskip

\noindent{\bf \large Acknowledgements}

\medskip

This work was supported by the National Key Research and Development Program of China (Grant No. 2016YFB0201304) and the Major Research Plan of National Natural Science Foundation of China (Grant Nos. 91430215, 91530323).

\bigskip
\bigskip

\noindent{\bf \large References}

\medskip

\small
{
	
[1] A.J. Hoffman, H.W. Wielandt, The variation of the spectrum of a normal matrix, Duke Math. J. 20 (1953) 37--39.

[2] J.-G. Sun, On the variation of the spectrum of a normal matrix, Linear Algebra Appl. 246 (1996) 215--223.

[3] Y. Song, A note on the variation of the spectrum of an arbitrary matrix, Linear Algebra Appl. 342 (2002) 41--46.

[4] J.-G. Sun, On the perturbation of the eigenvalues of a normal matrix, Math. Numer. Sinica 3 (1984) 334--336.

[5] Z. Zhang, On the perturbation of the eigenvalues of a non-defective matrix, Math. Numer. Sinica 6 (1986) 106--108.

[6] W. Li, W. Sun, The perturbation bounds for eigenvalues of normal matrices, Numer. Linear Algebra Appl. 12 (2005) 89--94.

[7] W.M. Kahan, Spectra of nearly Hermitian matrices, Proc. Amer. Math. Soc. 48 (1975) 11--17.

[8] P. Henrici, Bounds for iterates, inverses, spectral variation and fields of values of non-normal matrices, Numer. Math. 4 (1962) 24--40.

[9] J.-G. Sun, On the Wielandt--Hoffman theorem, Math. Numer. Sinica 2 (1983) 208--212.

[10] W. Li, S.-W. Vong, On the variation of the spectrum of a Hermitian matrix, Appl. Math. Lett. 65 (2017) 70--76.

[11] R.A. Horn, C.R. Johnson, Topics in Matrix Analysis, Cambridge University Press, Cambridge, 1991.

[12] S.C. Eisenstat, I.C.F. Ipsen, Three absolute perturbation bounds for matrix eigenvalues imply relative bounds, SIAM J. Matrix Anal. Appl. 20 (1998) 149--158.

[13] I.C.F. Ipsen, Relative perturbation results for matrix eigenvalues and singular values, Acta Numer. 7 (1998) 151--201.

[14] I.C.F. Ipsen, A note on unifying absolute and relative perturbation bounds, Linear Algebra Appl. 358 (2003) 239--253.

[15] R.-C. Li, Relative perturbation theory: I. Eigenvalue and singular value variations, SIAM J. Matrix Anal. Appl. 19 (1998) 956--982.

}

\end{document}